\documentclass[a4paper,12pt]{article}
\usepackage{epsfig}
\usepackage{amsmath}
\usepackage{amssymb}
\usepackage{mathtools}
\usepackage[top=2.5cm, bottom=2.5cm, left=1.5cm, right=1.5cm]{geometry}
\usepackage{microtype}
\usepackage{lmodern}
\usepackage[utf8]{inputenc}
\usepackage{float}
\usepackage{bera}
\DeclarePairedDelimiter{\ceil}{\lceil}{\rceil}

\newcommand{\ds}{\displaystyle }

\newcommand{\beq}{\begin{equation} }
\newcommand{\eeq}{\end{equation}}
\author{{\Large Thomas M. Michelitsch$^a$, Federico Polito$^b$,
Alejandro P. Riascos$^c$ }\\ \\ 
\footnotesize{$^a$ Sorbonne Universit\'e, Institut Jean le Rond d’Alembert, 
CNRS UMR 7190} \\
\footnotesize{4 place Jussieu, 75252 Paris cedex 05, France} \\ 
\footnotesize{E-mail: michel@lmm.jussieu.fr}\\[1ex]
\footnotesize{$^b$ Department of Mathematics ``Giuseppe Peano'', University of Torino, Italy}  \\ 
\footnotesize{E-mail: federico.polito@unito.it}\\[1ex]
\footnotesize{$^c$ Instituto de F\'isica, Universidad Nacional Aut\'onoma de M\'exico} \\[1ex]
\footnotesize{Apartado Postal 20-364, 01000 Ciudad de M\'exico, M\'exico}  \\
\footnotesize{E-mail: aperezr@fisica.unam.mx} \\[4ex]
}
\title{Semi-Markovian discrete-time telegraph process with generalized Sibuya waiting times}
\begin{document}

\maketitle
\begin{abstract}
In a recent work we introduced a semi-Markovian discrete-time generalization of the telegraph process. We referred this random walk to as `{\it squirrel random walk}' (SRW).
The SRW is a discrete-time random walk on the one-dimensional infinite lattice where the
step direction is reversed at arrival times of a discrete-time renewal process
and remains unchanged at uneventful time instants.
We first recall general notions of the SRW. 
The main subject of the paper is the study of the SRW where the step direction switches at the arrival times of a generalization of the Sibuya discrete-time renewal process (GSP) which only recently appeared in the literature. The waiting time density of the GSP, the `generalized Sibuya distribution' (GSD) is such that the moments are finite up to a certain order $r\leq m-1$ ($m \geq 1$) and diverging
for orders $r \geq m$ capturing all behaviors from broad to narrow and containing the standard Sibuya distribution as a special case ($m=1$). We also derive some new representations for the generating functions related to the GSD.
We show that the generalized Sibuya SRW exhibits several regimes of anomalous diffusion depending on the lowest order $m$ of diverging GSD moment.
The generalized Sibuya SRW opens various new directions in anomalous physics.
 \\[1mm]
\end{abstract}
\newpage
\tableofcontents
\section{Introduction}
\label{Introduction}
The telegraph process is an important model for transport where the velocity of the moving particle remains finite with a wide range of applications in physically existing and observable transport phenomena \cite{GionaCairoliKlages2020}.
The classical telegraph process (Poisson-Kac process) \cite{Goldstein1951,Kac1974} is defined as a 
one-dimensional motion of a particle
with constant velocity where the velocity direction is switched randomly at Poisson renewal times.
The classical telegraph process is Markovian inheriting this feature from the Poisson process. Its time evolution is governed by the (hyperbolic) telegrapher's (also called Cattaneo) equation avoiding physically forbidden infinite propagation velocities of the moving particle as occurring in the parabolic
standard diffusion equation.
Meanwhile, a wide range of semi-Markovian variants of the telegraph process, including fractional generalizations were developed to model anomalous transport (see among others \cite{mirko2,mirko,ga,Gorskaetal2020,Masoliver2016,orsbeg}) as well as
a tempered space-fractional generalization \cite{Beghin-etal2022}.
Compte and Metzler considered phenomenological fractional generalizations 
\cite{CompteMetzler1997,CompteMetzler1999} and related this model with the Montroll--Weiss 
continuous-time random walk (CTRW) framework. They found ballistic behavior for long times when the waiting time distribution has diverging mean, and enhanced non-ballistic transport in cases in which the waiting time distribution has a finite mean. Such behavior also occurs in the large-time asymptotics in our recent SRW model \cite{SRW2022} and in the model studied in the present paper. 
Further works considered the occurrence of random velocities \cite{StadjeZacks2004}, a relativistic model and analysis of occupation times, respectively \cite{BeghinNiedduOrsingher2001,BogachevRatanov2011},
Erlang distributed velocity reversals \cite{DiCrescendo2001,CrescenzoMartinucci2013}, distribution of the maximum \cite{cinque} --- consult also the references therein.

These works refer to continuous-time variants of the telegraph process. On the other hand many real-world datasets for instance in finance refer to discrete observation times \cite{deGregorioIacus2008}.
Therefore, it appears natural to consider discrete-time variants of (generalized) telegraph type processes calling still for thorough analytical investigation. 
In a recent work \cite{SRW2022}, we introduced a discrete-time semi-Markovian version of the telegraph process, the `{\it squirrel random walk}' (SRW), which is also subject of the present paper.
We chose that name since the SRW walker (the `squirrel') in a sense has a `weaker' memory as in walks with a full memory of their history such as the `elephant' walker in the
so called elephant random walk (ERW) 
\cite{SchuetzTrimper2004}. 

Our paper is organized as follows. In Section \ref{prelim} we give a brief account for discrete-time renewal processes and introduce pertinent generating functions for the present study. These generating functions
will be used in Section \ref{SRW} where we give an outline of basic notions of the SRW. Section \ref{general_Sibuya} is devoted to a generalization of the Sibuya distribution which appeared in the literature only recently \cite{KozubowskiPodgorski2018}. The speciality of the `generalized Sibuya distribution' (GSD) is that it has existing integer order moments only up to a certain order. In this way the GSD covers a wide range 
of behaviors from narrow to broad.
In Section \ref{Gen_Sib_SRW} we analyze the SRW where the step direction is
reversed at generalized Sibuya arrival times. We call this walk the `generalized Sibuya SRW'. The (anomalous) diffusive features of this walk are analyzed in Section \ref{anomalous_transport}.
\section{Discrete-time renewal process and related generating functions}
\label{prelim}
The theory of renewal processes generally provides a simple but powerful framework within the theory of point processes in the case in which the waiting times between the events are independent. It is well known that for the Markovian 
(memoryless) cases the waiting times can be either exponentially distributed when time is continuous or geometrically distributed for discrete time. However, in many real world situations the Markov property does no longer hold and the process has a memory. This leads to the study of semi-Markov processes \cite{Cox1962,Feller1993,Levy1954,Pyke1961,Smith1955} (and see the references therein).
\\
Most of the mentioned models, and many others, deal with renewal processes defined for continuous time. Compared to their continuous-time counterparts, the appearance of discrete-time renewal processes in the literature is relatively rare. Based on the theory of semi-Markov processes
Barbu and Limnios analyzed \cite{BarbuLimnios2008} discrete-time renewal processes and
recent works \cite{MichelitschPolitoRiascos2021,ADTRW2021,PachonPolitoRicciuti2021} are devoted to the analysis and applications of discrete-time renewal processes. In the present section we recall some of their features as far we need them for the SRW model.
We consider a discrete-time counting (renewal) process as follows:
\beq
\label{discrete-time}
{\cal N}(t) = \max(n\geq 0: J_n \leq t)  , \hspace{1cm} {\cal N}(0)=0 , \hspace{1cm} t =0,1,2\ldots
\eeq
The arrival times (renewal times) $J_n \in \mathbb{N} = \{1,2,\ldots\}$ (time instants of events, arrivals) are characterized by the random variables
\beq
\label{renewal_chain}
J_n=\sum_{j=1}^n \Delta t_j , \hspace{0.5cm} J_0=0, \hspace{0.5cm} \Delta t_j  \in \mathbb{N}
\eeq
with IID (independent and identically distributed) strictly positive integer increments $\Delta t_j \geq 1 $ (`interarrival times'  or `waiting times' in the renewal interpretation). The renewal chain
(\ref{renewal_chain}) is a discrete version of a strictly increasing subordinator. We refer to the recent article \cite{PachonPolitoRicciuti2021} elaborating essential elements of the related theory of discrete-time semi-Markov processes. 
The increments follow a discrete-time probability density function (PDF)
\beq
\label{first_arrival}
\mathbb{P}(\Delta t=k) = \psi_k , \hspace{1cm}  k=1,2,\ldots
\eeq
supported on positive integers $k\in \mathbb{N}$ with $\psi_0=0$ ensuring strictly positive waiting times. We employ the terms PDF and `density' in both cases, discrete-time and continuous-time.
The inverse of the renewal chain (\ref{renewal_chain}) is
the discrete-time counting process (\ref{discrete-time}) which counts the events (renewals) up to time $t$.
In the present paper we extensively use generating functions (GFs).
It is useful to introduce the GF of the waiting time density
\beq
\label{gen_fu}
\left\langle u^{\Delta t} \right\rangle= {\bar \psi}(u) = \sum_{t=1}^{\infty}\psi_tu^t , \hspace{1cm} |u| \leq 1,
\eeq
which fulfills ${\bar \psi}(u)\big|_{u=1}=1$ indicating normalization of (\ref{first_arrival}) and be reminded that $\psi(t)$ is supported on non-zero integers $t\in \mathbb{N}$.
Generally, the notation 
\beq
\label{GF_gen_notation}
 {\bar f}(u) = \sum_{t=0}^{\infty} f(t)u^t
 \eeq
(with suitably chosen $u$) stands
for the generating function (GF) of discrete functions $f(t)$ supported on
$t\in \mathbb{N}_0$.
We employ $\mathbb{E}(\ldots) =\langle (\ldots)\rangle$ as equivalent notations for expectation values
where we will often use
\beq
\label{average_f}
\left\langle f(\Delta t) \right\rangle = \sum_{r=1}^{\infty}f(r) \psi_r .
\eeq
Convenient is to introduce the indicator function 
\beq
\label{indicator_function}
\Theta(J_n,t,J_{n+1}) = \left\{\begin{array}{clr} \ds 1 & \ds {\rm if} & \hspace{0.25cm}\ds J_n\leq t \leq J_{n+1}-1,\hspace{0.25cm}
\\ \\
\ds 0  & \ds  {\rm otherwise} &
 \end{array}\right.
\eeq
which is one for ${\cal N}(t)=n$ and null else.
Then, the `state probabilities' 
(probabilities for $n$ arrivals up to time $t$) are given by \cite{SRW2022} (and see \cite{GodrecheLuck2001} for a related analysis for continuous time renewal processes)
\beq
\label{state_prob}
\mathbb{P}[{\cal N}(t)=n] = \phi^{(n)}(t) =  \left\langle\Theta(J_n,t,J_{n+1})\right\rangle.
\eeq
A quantity of interest is the variable $B_{n,t}=t-J_n$ containing information on the persistence of ${\cal N}(t)$ in state $n$ and which gives a connection to the
`aged renewal process' \cite{SRW2022} (and consult \cite{BarkaiCheng2003,GodrecheLuck2001,SchulzBarkaiMetzler2014} for the continuous time cases). We have (read $\mathbb{P}(A_1|A_2)$ as the probability of $A_1$ conditional to $A_2$)
\beq
\label{cond_prob}
f_B(\tau,t,n) = \mathbb{P}[B_{n,t}=\tau|{\cal N}(t)=n] = \left\langle \delta_{\tau,t-J_n} \Theta(J_n,t,J_{n+1})\right\rangle .
\eeq
We consider its double GF ${\bar f}_B(w,u,n)= \sum_{\tau=0}^{\infty}\sum_{t=0}^{\infty} f_B(\tau,t,n) w^{\tau}u^t$ ($|u|<1, |w|\leq 1$) which yields
\beq
\label{d_GFf_B}
{\bar f}_B(w,u,n) = \left\langle w^{-J_n} \sum_{t=J_n}^{J_{n+1}-1} (wu)^t   \right\rangle =  \left\langle u^{\Delta t_1+\ldots+ \Delta t_n}\frac{1-(uw)^{\Delta t_{n+1}}}{1-uw} \right\rangle =
[{\bar \psi}(u)]^n\frac{1-\bar \psi(uw)}{1-uw}
\eeq
where ${\bar f}_B(1,u,n) = [{\bar \psi}(u)]^n\frac{1-\bar \psi(u)}{1-u}$ recovers the GF of the state probabilities. In these derivations
we always use the IID feature of the $\Delta t_j$ and (\ref{gen_fu}).
The following two relations 
are related with the SRW propagator, namely
\beq
\label{average_cases}
g(t;\zeta_1,\ldots,\zeta_n;\zeta_{n+1})= \left\langle \zeta_{n+1}^{t-J_n+1}\Theta(J_n,t,J_{n+1}) \zeta_1^{\Delta t_1 -1}   \prod_{j=2}^n\zeta_j^{\Delta t_j}\right\rangle , \hspace{1cm} n,t \in \mathbb{N}_0 , \hspace{0.5cm} |\zeta_j| \leq 1
\eeq
recovering for $\zeta_j=1$ the state probabilities and 
\beq
\label{gt_all_zetas}
g(t,\{\zeta_j\}) = \sum_{n=0}^{\infty}g(t;\zeta_1,\ldots,\zeta_n;\zeta_{n+1}).
\eeq
The function $g(t;\zeta_1,\ldots,\zeta_n;\zeta_{n+1})$ in equation (\ref{average_cases}) has the GF
\beq
\label{GF_mean}
\begin{array}{clr}
\ds {\bar g}(u;\zeta_1,\ldots,\zeta_n;\zeta_{n+1}) & = \ds \sum_{t=0}^{\infty} u^t f(t;\zeta_1,\ldots,\zeta_n;\zeta_{n+1}) = \zeta_1^{-1}\zeta_{n+1} \left\langle \prod_{j=1}^n\zeta_j^{\Delta t_j}\sum_{t=J_n}^{J_{n+1}-1} \zeta_{n+1}^{t-J_n} u^t \right\rangle & \\ \\
 & = \ds \zeta_1^{-1}\zeta_{n+1}\frac{1-{\bar \psi}(u\zeta_{n+1})}{1-u\zeta_{n+1}} \prod_{j=1}^n {\bar \psi}(u\zeta_j) &
 \end{array}  (|u|<1)
\eeq
with ${\bar g}(u,\{\zeta_j\})\big|_{\zeta_j=1,u=1} =\frac{1}{1-u}$ corresponding to
the normalization of the state probabilities 
$$g(t,\{\zeta_j\})\big|_{\zeta_j=1} = \sum_{n=0}^{\infty}\phi^{(n)}(t) =1.$$
For the SRW the particular case when $\zeta_j$ alternate as $\zeta_{2j+1}=\zeta_1$ and $\zeta_{2j}= \zeta_2$ is pertinent.
Then we have for (\ref{GF_mean}) the GF
\beq
\label{special-SRW-case}
\ds {\bar g}_n(u;\zeta_1,\zeta_2)= \left\{\begin{array}{llr} 
\ds \frac{1-{\bar \psi}(u\zeta_1)}{1-u\zeta_1}\left[{\bar \psi}(u\zeta_1){\bar \psi}(u\zeta_2)\right]^{\ell}  , & n=2\ell & \\ \\ \ds \zeta_1^{-1}\zeta_2
\frac{1-{\bar \psi}(u\zeta_2)}{1-u\zeta_2}
{\bar \psi}(u\zeta_1)\left[{\bar \psi}(u\zeta_1){\bar \psi}(u\zeta_2)\right]^{\ell} , & n=2\ell+1. &
\end{array}\right.   (\ell=0,1,2,\ldots).
\eeq
Summation over $n$ yields the GF of (\ref{gt_all_zetas})
as
\beq
\label{GF_of_prop}
{\bar g}(u;\zeta_1,\zeta_2) = \left[\frac{1-{\bar \psi}(u\zeta_1)}{1-u\zeta_1} +\zeta_1^{-1}\zeta_2 {\bar \psi}(u\zeta_1) \frac{1-{\bar \psi}(u\zeta_2)}{1-u\zeta_2}\right]\, \frac{1}{1-{\bar \psi}(u\zeta_1){\bar \psi}(u\zeta_2)}.
\eeq
We will come back to these GFs in the context of the SRW propagator in the subsequent section.

\section{The squirrel random walk -- SRW}
\label{SRW}
Here we give a brief outline of the `{\it squirrel random walk}' (SRW), for an extensive study we refer to our recent work \cite{SRW2022}.
The SRW is a discrete-time random walk $X_t \in \mathbb{Z}$ where directed unit steps $\sigma_t \in \{-1,1\}$ are performed at integer time instants (we denote
with $t \in \mathbb{N}_0$ the time variable) 
\beq
\label{random_var}
X_t = \sum_{r=1}^t \sigma_r , \hspace{.5cm} t=1,2,\ldots, \hspace{1cm} X_0=0.
\eeq
The directions of steps are switched at arrival times of a discrete-time renewal process ${\cal N}(t)$.
A precise definition of the SRW is as follows:
\begin{enumerate}
\item [(i)] At uneventful time instants $t$,
the squirrel performs a unit step $\sigma_t=\sigma_{t-1}$ in the same direction as at $t-1$ where this holds for $t\geq 2$.
\item [(ii)] At arrival times $t$, 
the squirrel changes the step direction with respect to the previous step $\sigma_t=-\sigma_{t-1}$. 
\item [(iii)] 
We define that no step is performed at $t=0$ in order to ensure the initial condition $X_0=0$. The first step is performed at $t=1$ in the direction
$\sigma_1 = {\tilde \sigma}_0$ if $t=1$ is uneventful and 
$\sigma_1 = - {\tilde \sigma}_0$ if there is an event at $t=1$.
The direction ${\tilde \sigma}_0$ can be thought as
either prescribed or randomly chosen.
\end{enumerate}
In the following, we consider ${\tilde \sigma}_0 \in \{-1,1\}$ as given.
From the above it follows that the steps can be represented as
\beq
\label{increment}
\sigma_t = {\tilde \sigma}_0[(-1)^{{\cal N}(t)} - \delta_{t0}] , \hspace{1cm} t \in \mathbb{N}_0
\eeq
where the Kronecker-$\delta_{t0}$ ensures that no step is performed
at $t=0$.
Therefore, given ${\cal N}(t) = n$,
\beq
\label{gen_SRW}
\begin{array}{clr}
\ds X_t   & =  \ds {\tilde \sigma}_0\left[-1+\Delta t_1 -\Delta t_2 +\ldots + 
 (-1)^{n-1}\Delta t_n + (-1)^n(t-J_n+1)\right]  & \\ \\
  & =  X_t^{(+)}-X_t^{(-)}  & 
  \end{array} 
\eeq
with initial condition $X_0=0$ and where 
$X_t^{(+)}, X_t^{(-)}$ cover the steps in ${\tilde \sigma}_0$- and in the opposite direction, respectively.
Now we introduce the propagator (probability that the squirrel at time $t$ is sitting on $X\in \mathbb{Z}$) as follows
\beq
\label{SRW_propagator}
\mathbb{P}[X_t=X] = P(X,t) = \left\langle \delta_{X,X_t} \right\rangle ,\hspace{1cm} X \in \mathbb{Z}, \hspace{0.5cm} t \in \mathbb{N}_0
\eeq
with the Kronecker symbol $\delta_{A,B}$. Now using
$$\delta_{A,B} = \frac{1}{2\pi}\int_{-\pi}^{\pi} e^{i\varphi(A-B)}{\rm d}\varphi, \hspace{1cm} A,B \in \mathbb{Z}$$ 
we have
\beq
\label{rewrite_SRW_propagator}
P(X,t) = \frac{1}{2\pi}\int_{-\pi}^{\pi} e^{i\varphi X} 
\left\langle e^{-i\varphi X_t }\right\rangle {\rm d}\varphi
\eeq
and with (\ref{gen_SRW}) and (\ref{GF_of_prop}) the characteristic function writes
\beq
\label{Fourier_propagator}
P_{\varphi}(t)  =
\left\langle e^{-i\varphi X_t }\right\rangle = 
g(t;e^{-i\varphi{\tilde \sigma}_0} ,e^{i\varphi{\tilde \sigma}_0})  , \hspace{1cm} \varphi \in [-\pi,\pi].
\eeq
In addition, the GF ${\bar P}_{\varphi}(u)= {\bar g}(u;e^{-i\varphi{\tilde \sigma}_0} ,e^{i\varphi{\tilde \sigma}_0})$ is useful where ${\bar P}_{\varphi}(u)\big|_{\varphi=0}=\frac{1}{1-u}$ tells us that the propagator $P(X,t)$ is a (spatially) normalized PDF.
Then we introduce 
\beq
\label{state_plynomial}
{\cal P}(v,t) = \langle v^{{\cal N}(t)} \rangle = \sum_{n=0}^t\mathbb{P}[{\cal N}(t)=n)]v^n
\eeq
which is a polynomial of degree $t$ (`state polynomial') since
$\mathbb{P}[{\cal N}(t)=n)] =0 $ for $n>t$ as ${\cal N}(t) \leq t$
with initial condition $\mathbb{P}[{\cal N}(t)=n)] = \delta_{n0}$.
The feature ${\cal P}(1,t)=1$ reflects normalization of the state probabilities, and for $v=-1$ the average step is contained, namely
\beq
\label{mean_sigma}
\langle \sigma_t \rangle = {\tilde \sigma}_0[ \left\langle (-1)^{{\cal N}(t)} \right\rangle -\delta_{t0}] = {\tilde \sigma}_0[{\cal P}(-1,t) - \delta_{t0}]
\eeq
where $\delta_{t0}$ takes into account that no step is performed at $t=0$
maintaining the initial condition $X_0=0$.
The GF of the average steps then takes
\beq
\label{mean_step_gen}
{\bar \sigma}(u) = \sum_{t=0}^{\infty} \langle \sigma_t\rangle u^t = 
{\tilde \sigma_0}\left[{\bar {\cal P}}(-1,u) - 1 \right]
\eeq
with the GF of the state polynomial
\beq
\label{state_gen}
{\bar {\cal P}}(v,u) = 
\frac{1-{\bar \psi}(u)}{(1-u)[1-v{\bar \psi}(u)]}
,\hspace{1cm} |u|< 1, \hspace{1cm} |v| \leq 1.
\eeq
The GF of the expected position $\langle X_t\rangle$ then reads
\beq
\label{gen_fu_of_position}
{\bar X}^{(1)}(u) = i \frac{\partial }{\partial \varphi}{\bar P}_{\varphi}(u)\big|_{\varphi=0}  =
\frac{{\bar \sigma}(u)}{1-u} =
 \frac{[1-{\bar \psi}(u)]{\tilde \sigma}_0 }{(1-u)^2[1+{\bar \psi}(u)]}-\frac{{\tilde \sigma}_0}{1-u}
\eeq
with the initial condition ${\bar X}^{(1)}(u)\big|_{u=0} = \langle X_0\rangle = 0$. We will focus on the `generalized Sibuya SRW' where the instants of the step reversals are drawn from a generalization of the Sibuya
distribution which is subject of the subsequent section.
\section{Generalized Sibuya counting process}
\label{general_Sibuya}
Here we consider a discrete 
time counting process ${\cal N}_{\lambda}(t)$ with IID generalized Sibuya waiting-times.
The resulting generalized Sibuya distribution (GSD) was to our knowledge first introduced and thoroughly studied by Kozubowski and Podg\'orski \cite{KozubowskiPodgorski2018}. 
In the present section we recall the GSD in the light of discrete-time renewal processes and derive some auxiliary GFs needed in Sections \ref{Gen_Sib_SRW} and \ref{anomalous_transport}.
To this end we introduce the GSP waiting time PDF as follows
\beq
\label{PDF_waiting}
\begin{array}{clr}
\ds \psi_{\lambda}(t) &=  \ds \frac{\Gamma(1-\lambda)\Gamma(m)}{\Gamma(m-\lambda)}(-1)^{t+m}  \left(\begin{array}{l} \,\,\, \lambda \\ 
t+m-1 \end{array}\right) &  \\ \\
 & = \ds \frac{\lambda \Gamma(m)}{\Gamma(m-\lambda)}\frac{\Gamma(m-\lambda+t-1)}{\Gamma(m+t)} &
 \end{array} \hspace{1cm} (t \in \mathbb{N}, \hspace{0.5cm} \lambda >0, \hspace{0.5cm} m=\ceil{\lambda})
\eeq
where $\ceil{\lambda}$ indicates the ceiling function, producing the smallest integer larger than or equal to $\lambda$
where we mainly focus on non-integer $\lambda$.
The positiveness of this expression is easily confirmed by accounting for $\lambda=m-1+\mu$, $\mu \in (0,1)$ (see especially (\ref{GSP_Pochhammer_rep})). For $\lambda \in (0,1)$, i.e. $m=1$ (\ref{PDF_waiting}) recovers the standard Sibuya distribution \cite{Sibuya1979}.
We refer the PDF (\ref{PDF_waiting}) as `{\it generalized Sibuya distribution}' (GSD).

We construct the `generalized Sibuya counting process' (GSP) such that it has a discrete
waiting-time PDF with finite moments up to order $r\leq m-1$ ($m\geq 1$) and diverging moments of orders $r\geq m$. 
This implies that
$\frac{d^r}{du^r}{\bar \psi}(u)\big|_{u=1}$ is finite for $r\leq m-1$ and $\frac{d^r}{du^r}{\bar \psi}(u)\big|_{u=1} \to \infty$ for $r \geq m$ ($m\geq 1$). 
The GSD has the following GF
\beq
\label{gen-Sib}
 {\bar \psi}_{\lambda}(u)  =  \frac{\Gamma(1-\lambda)\Gamma(m)}{\Gamma(m-\lambda)}u^{1-m}[H_{\lambda}(u)-(1-u)^{\lambda}] , \hspace{0.5cm} 0\leq m-1 < \lambda < m  \in \mathbb{N} , \hspace{0.5cm} |u| \leq 1
\eeq
with 
\beq
\label{Hlam}
H_{\lambda}(u) = \sum_{r=0}^{m-1} (-1)^r
\left(\begin{array}{l} \lambda \\ r \end{array}\right) u^r
\eeq
which removes the terms with alternating signs in the expansion $-(1-u)^{\lambda}$ thus ${\bar \psi}_{\lambda}(u)$ contains all non-alternating orders $u^r$ for $r\geq m$ of this expansion. One can easily verify that the sign of these terms is 
$(-1)^{m-1}=\mathrm{sign}(H_{\lambda}(1))$ thus
(\ref{gen-Sib}) contains only non-negative coefficients and is of the form 
$${\bar \psi}_{\lambda}(u) = u^{1-m} (H_{\lambda}(1))^{-1}[H_{\lambda}(u)-(1-u)^{\lambda}] = u^{1-m} {\bar g}_{\lambda}(u)$$ where $u^{1-m}$ shifts the distribution $g_{\lambda}(t)$ by $m-1$ to the left, ensuring that 
$\psi_{\lambda}(t)=g_{\lambda}(t+m-1)$ is nonzero from $t\geq 1$.
The normalization factor is obtained as
$$H_{\lambda}(1) =
\frac{1}{(m-1)!}\frac{d^{m-1}}{du^{m-1}}(1-u)^{\lambda-1}\big|_{u=0} =  (-1)^{m-1}\left(\begin{array}{l}  \lambda-1 \\ m-1 \end{array}\right) =\frac{\Gamma(m-\lambda)}{\Gamma(1-\lambda)\Gamma(m)} = \frac{(1-\lambda)_{m-1}}{(1)_{m-1}}.$$ 
Let us remark that, although $\lambda \notin \mathbb{N}$, integer values $\lambda = m$ are admissible
retrieving ${\bar \psi}_{m}(u) = u$ corresponding to the trivial (deterministic) counting process ${\cal N}_{m}(t)=t$ and
coinciding with the limit $p\to 1-$ of a Bernoulli counting process where $p$ indicates the probability of a Bernoulli success.
By construction ${\bar \psi}_{\lambda}(u) =O(u)$ and we employ the Pochhammer symbol
$$
(a)_k = a(a+1)\ldots (a+k-1) = \left\{\begin{array}{ll} \ds \frac{\Gamma(a+k)}{\Gamma(a)} , & k \in \mathbb{N} \\ \\
\ds  1 , & k=0
\end{array}
\right.
$$
and we mention the useful features $(a)_{k+r} = (a)_k (a+k)_r$.
We notice for later use the property
\beq
\label{propHlam}
\frac{d^{\ell}}{du^{\ell}}H_{\lambda}(u) = (-\lambda)_{\ell}  H_{\lambda-\ell}(u) 
\eeq
thus $\frac{1}{H_{\lambda}(1)}\frac{d}{du}H_{\lambda}(u)\big|_{u=1}= \frac{\lambda(m-1)}{\lambda-1}$ ($\lambda>1$) and $H_{\lambda-\ell}(u)=0$ for $\ell \geq m$. Then we can expand (\ref{Hlam}) with respect to
$u-1$ as follows
\beq
\label{Hlam_expan}
H_{\lambda}(u)= \sum_{\ell=0}^{m-1} \frac{(u-1)^{\ell}}{\ell !}\frac{d^{\ell}}{du^{\ell}}H_{\lambda}(u)\big|_{u=1} = \Gamma(m-\lambda) \sum_{\ell=0}^{m-1} \left(\begin{array}{l} \lambda \\ \ell \end{array}\right) 
\frac{(1-u)^{\ell}}{\Gamma(1+\ell-\lambda)\Gamma(m-\ell)}.
\eeq
Therefore, the GF (\ref{gen-Sib}) writes compactly as
\beq
\label{GSD_SibGF_expand}
{\bar \psi}_{\lambda}(u) =u^{1-m}{\bar g}_{\lambda}(u) =
u^{1-m}\left(\, 1-\frac{(1)_{m-1}}{(1-\lambda)_{m-1}} (1-u)^{\lambda} +
\sum_{\ell=1}^{m-1} \left(\begin{array}{l} \lambda \\ \ell \end{array}\right) 
\frac{(m-\ell)_{\ell}}{(1-\lambda)_{\ell}}(1-u)^{\ell} \,\right).
\eeq
In the large time limit we have for (\ref{PDF_waiting}) the representation (we employ symbol $\sim$ for asymptotic equality and use $\frac{\Gamma(t+a)}{\Gamma(t+b)} \sim t^{a-b}$ as $t\to \infty$)
\beq
\label{large_time_pdf}
\psi_{\lambda}(t) \sim \frac{\lambda \Gamma(m)}{\Gamma(m-\lambda)}\, t^{-\lambda-1} 
\eeq
which holds for any $\lambda>0$. Thus the GSD (\ref{PDF_waiting}) covers any power-law 
from narrow (large $\lambda$) to broad (small $\lambda$).
Especially, for $\lambda \in (0,1)$ ($m=1$) (\ref{PDF_waiting}) recovers the fat-tailed (broad) standard Sibuya distribution with diverging first moment \cite{PachonPolitoRicciuti2021,Sibuya1979}.
The long time asymptotics (\ref{large_time_pdf}), by invoking Tauberian theorems, is obtained from the asymptotic expansion of (\ref{GSD_SibGF_expand}) for $u\to 1-$ with the relevant part ${\bar \psi}_{\lambda}(u) \sim 1-\frac{(1-u)^{\lambda}}{H_{\lambda}(1)}$ (in which we can safely neglect the integer powers $(1-u)^n$ ($n>0$) as they do not have distributions with a long tail).

It appears instructive to consider the GSP in the light of a sequential trial scheme which can be adopted for any discrete-time renewal process (see \cite{ADTRW2021,PachonPolitoRicciuti2021} for details).
Perform a sequence of $k=1,2,\ldots \in \mathbb{N}$
(GSP-) trials where each trial has two possible outcomes, ‘‘success’’ or ‘‘fail’’ where we introduce the random variables 
$Z_k \in \{ 0, 1 \} $ ($k\geq 1$)
with $Z_k=1$ if the outcome is a success and $Z_k=0$ for a fail  and $Z_0=0$ (no trial at $t=0$). Then introduce the conditional probability  
$\alpha_k =\mathbb{P}[Z_k=1|\{Z_{r}=0\}_{r<k}]$ of a success in the $k$th
trial given there was no success in earlier trials.
Then performing at each integer time instant $t$ 
a trial we have for the GSP counting variable 
\beq
\label{GSP_countvar}
{\cal N}_{\lambda}(t) = \sum_{k=1}^tZ_k.
\eeq
Then the waiting time density $\psi_{\lambda}(t)$ has the interpretation as the probability of the first GSP-success 
at time $t$ \cite{ADTRW2021}, i.e.
\beq
\label{psi_lam_firstsucess}
\psi_{\lambda}(t) = \alpha_t \, (1-\alpha_1)\ldots(1-\alpha_{t-1}) =
\alpha_t \, S_{t-1}.
\eeq
Here $S_{k}= \prod_{r=1}^k (1-\alpha_r) = \sum_{r=k+1}^{\infty} \psi_{\lambda}(r)$ is the probability of a sequence of $k$ GSP fails (probability of no GSP success in $k$ trials, `survival probability'). 
We point out that any discrete PDF $\psi(t)$ indeed can be represented by such a sequential trial scheme with $\alpha_t = \psi(t)/[\sum_{r=t}^{\infty}\psi(r)$ (see \cite{PachonPolitoRicciuti2021} for details).

Let us elaborate this structure for the GSP. Unlike in a Bernoulli trial process
(characterized by the memoryless property $\alpha_t=p$ independent of $t$) the GSD has a memory which
is reflected by (\ref{psi_lam_firstsucess}) containing the complete history
up to this first GSP success. 
Then we can rewrite (\ref{PDF_waiting}) in terms of Pochhammer symbols as follows
\beq
\label{GSP_Pochhammer_represent}
\psi_{\lambda}(t)= \lambda \frac{(m-\lambda)_{t-1}}{(m)_t}.
\eeq 
Now, since $(a)_k= (a)_{k-1} (a+k-1)$, we have ($m=\ceil{\lambda}$)
\beq
\label{GSP_Pochhammer_rep}
\psi_{\lambda}(t)= \frac{\lambda}{m+t-1} \frac{(m-\lambda)_{t-1}}{(m)_{t-1}} = \frac{\lambda}{m+t-1}\left(1-\frac{\lambda}{m}\right)
\left(1-\frac{\lambda}{m+1}\right)\ldots \left(1-\frac{\lambda}{m+t-2}\right)
\eeq 
coinciding with the representation which ad hoc was introduced by Kozubowski and Podg\'orski \cite{KozubowskiPodgorski2018} and has clearly the structure (\ref{psi_lam_firstsucess}) where we identify
\beq
\label{structure_trialscheme}
\alpha_k = \frac{\lambda}{m+k-1} , \hspace{1cm} m, k\geq 1
\eeq
with $\alpha_1= \psi_{\lambda}(t)\big|_{t=1} = \frac{\lambda}{m}$.
Indeed $m=1$ ($\lambda=\mu \in (0,1)$) retrieves the standard Sibuya process.
We notice that if $\lambda = m \in \mathbb{N}$ the trivial counting process ${\cal N}_m(t)=t$ is recovered where each trial is a success with $\psi_{m}(t) = \delta_{t1}$.
From (\ref{GSP_Pochhammer_rep}) we obtain the `survival probability', i.e.\ the probability of no event up to time $t$ in a GSP (which we now denote with $\mathbb{P}[{\cal N}_{\lambda}(t)=0]= S_t = \phi_{\lambda}^{(0)}(t)$) 

\beq
\label{GSP_survival}
\phi_{\lambda}^{(0)}(t) = \ds \prod_{r=0}^{t-1}\left(1-\frac{\lambda}{m+r}\right) = \frac{(m-\lambda)_t}{(m)_t}  = \ds  \frac{\Gamma(m)}{\Gamma(m-\lambda)}\frac{\Gamma(m-\lambda+t)}{\Gamma(m+t)}  
\eeq
with initial condition $\phi_{\lambda}^{(0)}(t)\big|_{t=0}=1$.
The large time asymptotics
is obtained as
\beq
\label{large_times}
\phi_{\lambda}^{(0)}(t) \sim \frac{\Gamma(m)}{\Gamma(m-\lambda)} t^{-\lambda} , \hspace{1cm} (t \to \infty)
\eeq 
where for $m=1$ and $\lambda=\mu \in (0,1)$ these relations recover the case of the standard Sibuya process. 

For what follows we recall the Gauss hypergeometric function defined as \cite{WhittakerWatson1927}
\beq
\label{hyper_geo_function}
_2F_1(a, b; c ; u) = \sum_{r=0}^{\infty} 
\frac{(a)_r (b)_r}{r! \, (c)_r}u^r , \hspace{1cm} a,b,c \in \mathbb{R}, \hspace{0.5cm} c \notin \mathbb{Z}_{\leq 0}.
\eeq
It is sufficient to consider here
$a,b,c >0$ with $a+b-c<0$ where for large $r$ the coefficients decay with a power-law as
$\frac{(a)_r (b)_r}{(c)_r r!} \sim const \, r^{a+b-c-1}$ ($r\to \infty$). In this case (\ref{hyper_geo_function}) converges for
$|u|\leq 1$ and the Gauss summation theorem holds \cite{WhittakerWatson1927}:
\beq
\label{Gauss_summation}
_2F_1(a, b; c ; 1) = \frac{\Gamma(c)\Gamma(c-a-b)}{\Gamma(c-a)\Gamma(c-b)} , \hspace{1cm} ( c > a+b).
\eeq
The GF of the GSP survival probability (\ref{GSP_survival})
has then the form
\beq
\label{survival_GSP_hyper_geo}
{\bar \phi}_{\lambda}^{(0)}(u) =\sum_{r=0}^{\infty} \frac{(1)_r}{r!}  \frac{(m-\lambda)_r}{(m)_r} u^r  =\,\, _2F_1(1, m-\lambda; m; u)
\eeq 
converging for $|u| <1$ for standard Sibuya $m=1$ and for 
$|u|\leq 1$ for $m>1$ (see the asymptotic relation (\ref{large_times})).
Using (\ref{GSP_survival}) with $\phi^{(0)}_{\lambda}(k-1)= \sum_{r=k}^{\infty}\psi_{\lambda}(r)$ and by using
(\ref{Gauss_summation}) we reconfirm (\ref{structure_trialscheme}), namely
\beq
\label{reconfirm_alpha_k}
\alpha_k =\frac{\psi_{\lambda}(k)}{\sum_{r=k}^{\infty}\psi_{\lambda}(r)} =\big( \,_2F_1(1,k+m-1-\lambda; k+m ; 1) \big)^{-1} =\frac{\lambda}{k+m-1}.
\eeq
We evaluate the first moment of the random variable $T$ having GSD  (\ref{PDF_waiting})
existing for $m\geq 2$ ($\lambda>1$),
which yields
\begin{figure}[t]
\centerline{\includegraphics[width=0.7\textwidth]{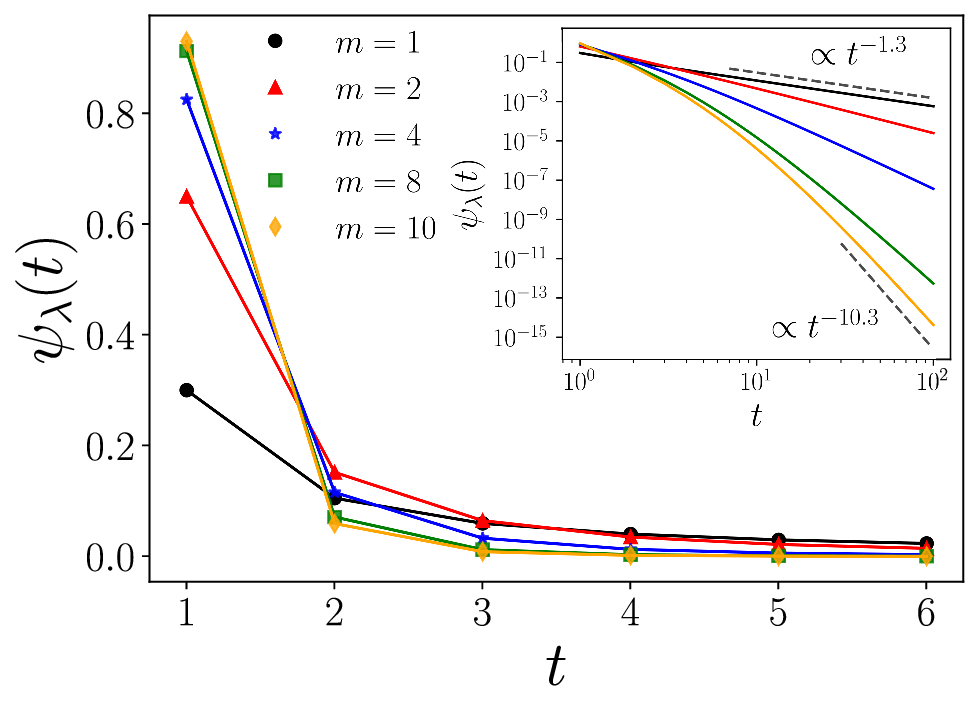}}
\vspace{-4mm}
\caption{\label{Fig1} Generalized Sibuya distribution $\psi_{\lambda}(t)$. We depict the GSD from Eq. (\ref{PDF_waiting}) for $\mu=0.3$ and different $m\geq 1$. The inset shows the results in logarithmic scale for $1\leq t\leq 100$, we present with dashed lines the power-law relation $\propto t^{-m-\mu}$ associated to the asymptotic result in Eq. (\ref{large_time_pdf})  for $m=1$ and $m=10$. }
\end{figure}
\beq
\label{first_moment}
\begin{array}{clr}
\ds \left\langle  T \right\rangle_{\lambda} &= \ds \frac{d}{du}{\bar \psi}_{\lambda}(u)\big|_{u=1} = \lambda \sum_{t=1}^{\infty} \frac{(m-\lambda)_{t-1}}{m_{t-1}} \frac{t}{t+m-1} & \\ \
& = \ds \lambda \sum_{t=1}^{\infty} \frac{(m-\lambda)_{t-1}}{m_{t-1}}(1-\frac{m-1}{m-1+t}) & \\ \\
& = \ds \lambda  {\bar \phi}^{(0)}(u)\big|_{u=1} -(m-1) & \\ \\
& = \ds \lambda \,\, _2F_1(1, m-\lambda; m; 1) -(m-1) = \frac{m-1}{\lambda-1} ,& \hspace{1cm} m \geq 2 ,\,\, (\lambda > 1)
\end{array}
\eeq
also conveniently obtained by accounting for representation (\ref{GSD_SibGF_expand}).
Note that, since $\psi_{\lambda}>0 $ on positive integers it is necessarily $\langle T \rangle_{\lambda} \ge 1$ which is fulfilled by this relation as $m=\ceil{\lambda} \ge \lambda$. We also observe that when we put $\lambda=m$ the first moment is consistent with the behavior of the 
corresponding trivial counting process (with $\psi(t)=\delta_{t1}$).
Further, for $\lambda \to \infty$ we have 
$$ \left\langle T \right\rangle_{\lambda} = \frac{1}{1+\frac{\mu-1}{m-1}} \to 1+$$ reflecting that, for large $\lambda$, the GSD becomes extremely narrow. Fig.~\ref{Fig1} shows the GSD for different values of $m$. Notice that with increasing $m$ the GSD becomes more narrow, and this is also reflected by the large time power-law  scaling (\ref{large_time_pdf}). 

\subsection{Bernoulli time-changed with GSP and scaling limits}
Before we return to the SRW
it appears instructive to highlight some connections of the GSP with pertinent counting processes which have recently appeared in the literature and to consider scaling limits to continuous-time. To this end we introduce the composed counting process ${\cal N}_B({\cal N}_{\lambda}(t))$ where ${\cal N}_B$ is a Bernoulli process and ${\cal N}_{\lambda}$ a GSP independent of ${\cal N}_B$. This composition is a Bernoulli counting process time-changed with a GSP where Bernoulli trials are performed at arrival times of the GSP (which describes a random clock).
For outlines on such compositions we refer to 
\cite{MichelitschPolitoRiascos2021,ADTRW2021,OrsingerPolito2012}.
The GF of the waiting time PDF of this composition is given by
\beq
\label{compos}
{\bar \chi}_{\lambda,\xi}(u) = {\bar \psi}_B[{\bar \psi}_{\lambda}(u)] = 
\frac{\xi {\bar \psi}_{\lambda}(u)}{\xi+1-{\bar \psi}_{\lambda}(u)}
\eeq
with the Bernoulli waiting time GF 
${\bar \psi}_B(u)=\frac{\xi u}{\xi+1-u}$ with $\xi=\frac{p}{1-p}$ where $p$ denotes the probability of success in each single Bernoulli trial.
The limit $p\to 1-$ (i.e. $\xi\to \infty$) 
${\cal N}_B({\cal N}_{\lambda}(t)) \to {\cal N}_{\lambda}(t)$ retrieves the GSP.
In the Sibuya case $m=1$ the composed process ${\cal N}_B({\cal N}_{\lambda}(t))$ contains the so called `fractional Bernoulli counting process' (of type B). The fractional Bernoulli counting process was introduced in \cite{PachonPolitoRicciuti2021} and has the waiting time GF 
\beq
\label{compos_fractber}
{\bar \chi}_{\lambda,\xi}(u) = {\bar \psi}_B[{\bar \psi}_{\lambda}(u)] =
\frac{\xi}{\xi +(1-u)^{\lambda}}[1-(1-u)^{\lambda}], \hspace{2cm} \lambda \in (0,1).
\eeq
Evoking Tauberian arguments, the long-time asymptotics of the waiting time density
of the composed process can be obtained by expanding (\ref{compos}) for $u \to 1-$ and considering only the lowest non-integer order in $1-u$ (see (\ref{GSD_SibGF_expand})), namely
\beq
\label{asymptotics}
{\bar \chi}_{\lambda,\xi}(u) = \frac{{\bar \psi}_{\lambda}(u)}{1+
\frac{1}{\xi}(1-{\bar \psi}_{\lambda}(u))} \sim  1-\frac{1}{p}(1-{\bar \psi}_{\lambda}(u)) \sim 1- \frac{1}{p H_{\lambda}(1)}(1-u)^{\lambda}
\eeq
where $1/p$ is the mean waiting time in a Bernoulli process. We skip in this asymptotic relation all integer orders in $1-u$. Hence, we get the density of the composed process 
\beq
\label{asymp_comp}
\chi_{\lambda,\xi}(t) = - \frac{(-\lambda)_t}{t! \, p H_{\lambda}(1)} \sim 
\frac{\lambda \Gamma(m)}{p\Gamma(m-\lambda)} t^{-\lambda-1} 
 , \hspace{1cm} (t \to \infty)
\eeq
having the same tail of the GSD (\ref{large_time_pdf})
up to the multiplier $1/p$ (Bernoulli mean).
Now we can define a well-scaled limit to continuous time $t \in h {\mathbb N}_0 \to \mathbb{R}_+$ where (see \cite{MichelitschPolitoRiascos2021} for a thorough outline of such continuum limit procedures) 
\beq
\label{Lplace_trafo}
{\hat \chi}_{\lambda,\xi_0}(s) = \lim_{h\to 0} 
{\bar \chi}_{\lambda,\xi(h)}(e^{-hs}) = 
\lim_{h\to 0} 
\frac{\xi_0 h^{\eta} {\bar \psi}_{\lambda}(e^{-hs})}{\xi_0 h^{\eta}+1-{\bar \psi}_{\lambda}(e^{-hs})}
\eeq
with the scaling assumption $\xi(h)=\xi_0h^{\eta}$ ($\xi_0>0$ is an arbitrary constant independent of the time increment $h$ and of the physical dimension $sec^{-\eta}$). The scaling exponent $\eta$ has to be chosen such that this limit exists. Then, accounting for
(\ref{GSD_SibGF_expand}) with $u=e^{-hs} \to 1-$ we have the asymptotic relation 
\beq
\label{lowest orders}
\ds {\bar \psi}_{\lambda}(e^{-hs}) \sim \left\{\begin{array}{l}\ds  1- h^{\lambda} s^{\lambda}+o(h^{\lambda}) , \hspace{2cm}  (\lambda \in (0,1)) \\ \\ \ds 
1- h \left\langle T \right\rangle_{\lambda} s +o(h)  , \hspace{2cm} (\lambda >1) \end{array}\right.
\eeq 
with $\left\langle T \right\rangle_{\lambda}$ given by (\ref{first_moment}).
Therefore, there exist only two possible limits for (\ref{Lplace_trafo}), namely for $m=1$ (standard Sibuya case) we have $\eta=\lambda \in (0,1)$
thus we obtain for that limit
\beq
\label{standard_ML}
{\hat \chi}_{\lambda,\xi_0}(s) = \frac{\xi_0}{\xi_0+s^{\lambda}} , \hspace{2cm} \lambda \in (0,1) 
\eeq
which is the Laplace transform of the Mittag-Leffler density, obtained in \cite{PachonPolitoRicciuti2021} as the continuous time limit of fractional Bernoulli 
to the fractional Poisson renewal process (see e.g.\ \cite{Laskin2003} among the many papers on the subject). 
For $m>1$ this limit exists only if we choose $\eta=1$ thus we get (by introducing the new constant time scale constant $\zeta_0=\xi_0/\left\langle T \right\rangle_{\lambda}$) 
\beq
\label{standard_exponential}
{\hat \chi}_{\lambda,\xi_0}(s) = \frac{\zeta_0}{\zeta_0+s}        ,\hspace{2cm}  \lambda > 1
\eeq
which is the Laplace transform of the exponential density 
$\chi_{\lambda,\xi_0}(t) = \zeta_0 e^{-\zeta_0 t}$ of the standard Poisson process. For $m>1$
the composition ${\cal N}_B({\cal N}_{\lambda}(t))$ converges in the above defined scaling limit to the standard Poisson process. 

These features also come into play when we consider the 
scaling limit of the expectation of the rescaled GSP renewal chain $n^{-\rho} J_n^{(\lambda)} \to J_{\lambda} $ for $n\to \infty$ 
(see Eq.\ (\ref{renewal_chain})) and choose exponent $\rho$ such that this limit exists
\beq
\label{scaling_subordinator}
\begin{array}{clr}
\ds \left\langle e^{-s J_{\lambda}} \right\rangle & = \ds \lim_{n\to \infty} \left\langle \exp{(-\frac{s J_n^{(\lambda)}}{n^{\rho}})} \right\rangle &   \\ \\  & = \ds \lim_{n\to \infty}  [{\bar \psi}_{\lambda}(e^{-\frac{s}{n^{\rho}}})]^n  & = \ds \left\{\begin{array}{ll} \ds  e^{- s^{\lambda}} & \ds \hspace{1cm} \lambda \in (0,1) \\ \\ \ds
e^{- s \left\langle T \right\rangle_{\lambda}} & \ds \hspace{1cm} \lambda >1  
\end{array}\right. 
\end{array}
\eeq
where we use the IID feature of the interarrival times 
together with (\ref{gen_fu}) and we
have to choose $\rho=\frac{1}{\lambda}$ for $m=1$ (standard Sibuya) and $\rho=1$ for $m>1$.
Hence $ n^{-\rho} J_n^{(\lambda)} \to  J_{\lambda} $ is a stable subordinator. 
\section{Generalized Sibuya SRW}
\label{Gen_Sib_SRW}
\subsection{Large time asymptotics of the expected squirrel position}
Here we explore the diffusive features of the generalized Sibuya SRW where the step directions are switched at GSP arrival times. To this end consider first 
the large time asymptotics of the expected position.
From (\ref{GSD_SibGF_expand}) we have in the asymptotic expansion three cases 
\beq
\label{utooneexpand}
{\bar \psi}_{\lambda}(u) = \left\{\begin{array}{ll} \ds 1-(1-u)^{\lambda} , & \hspace{0.3cm} 0 < \lambda < 1  \\ \\ \ds 
1 -\left\langle T \right\rangle_{\lambda}(1-u) -\frac{1}{H_{\lambda}(1)}(1-u)^{\lambda} + o_I[(1-u)]+o[(1-u)^{\lambda}] , & \hspace{0.3cm} 1 < \lambda  < 2 \\ \\
\ds 1 -\left\langle T \right\rangle_{\lambda}(1-u) +  B_2^{(\lambda)}(1-u)^2 -\frac{1}{H_{\lambda}(1)}(1-u)^{\lambda} +o[(1-u)^2] , & \hspace{0.3cm} \lambda >2
\end{array}\right. (u \to 1-).
\eeq
We denote with the symbol 
$o_I[(1-u)]= a_2(1-u)^2+a_3(1-u)^3..$ a power series in $(1-u)$ containing solely integer powers of orders larger than one.
Further we use $o[(1-u)^{\lambda}] \ll (1-u)^{\lambda}$ as $u \to 1-$.
The constant $B_2^{(\lambda)}$ is positive (existing for $\lambda >2$) and yields with (\ref{GSD_SibGF_expand})
\beq
\label{evalB}
\begin{array}{clr}
\ds B_2^{(\lambda)} &  = \ds  \frac{1}{2!}\frac{d^2}{du^2}{\bar \psi}_{\lambda}(u)\big|_{u=1}, \hspace{2cm} (\lambda > 2, m \geq 3) & \\ \\ 
& = \ds 
\frac{1}{2!}\left((1-m)(-m)u^{-m-1}{\bar g}_{\lambda}(u)+2(1-m) u^{-m}\frac{d}{du}{\bar g}_{\lambda}(u) + u^{1-m}\frac{d^2}{du^2}{\bar g}_{\lambda}(u)\right)\big|_{u=1} & \\ \\  
 & = \ds \frac{(m-1)(m-\lambda)}{(\lambda-1)(\lambda-2)} = \langle T \rangle_{\lambda} \frac{m-\lambda}{\lambda-2} ,  &
\end{array}
\eeq
where the non-negativeness of $B_2^{(\lambda)}$ can be seen from $m=\ceil{\lambda}$.
Since $B_2^{(\lambda)}= \sum_{t=1}^{\infty}\frac{t(t-1)}{2}\psi_{\lambda}(t)$ this coefficient contains also the second moment of the GSD 
\beq
\label{second_GSD_moment}
\langle T^2 \rangle_{\lambda} = 2B_2^{(\lambda)} +\langle T \rangle_{\lambda} = \frac{\langle T \rangle_{\lambda}(2m-\lambda-2)}{\lambda-2} = \frac{(m-1)(2m-\lambda-2)}{(\lambda-1)(\lambda-2)} , \hspace{1cm} (\lambda >2)
\eeq
where $\langle T \rangle_{\lambda}$ was determined in (\ref{first_moment}).
Further of interest is the variance existing for $\lambda >2$ ($m\geq 3$) which yields 
\beq
\label{variance_GSD}
{\cal V}_{\lambda} = \langle T^2 \rangle_{\lambda} - [\langle T \rangle_{\lambda}]^2 = 
\frac{\lambda(m-1)(m-\lambda)}{(\lambda-1)^2(\lambda-2)} = \frac{\lambda}{\lambda-1} B_2^{(\lambda)}.
\eeq
The GSP variance (\ref{variance_GSD}) coincides with the expression given in \cite{KozubowskiPodgorski2018} (see there Definition 1 with Remark 4 and Eq. (22) in that paper and identify $m=\nu+1$, $\lambda=\alpha$ in their notation).
Recall that we mainly consider $\lambda \notin \mathbb{N}$ and bear in mind that integer values $\lambda=m$ are admissible defining the deterministic counting process ${\cal N}_{m}(t) = t$ where for $\lambda\to m-$ ($m\geq 3$) the variance (\ref{variance_GSD}) exists and is vanishing.

Now with (\ref{gen_fu_of_position}) and (\ref{utooneexpand}) we obtain for the GF of the expected squirrel position
the asymptotic relation
\beq
\label{GS_SRW_position}
 {\bar X}_{\lambda}(u)  =  \left\{\begin {array}{ll}
 \ds  \frac{{\tilde \sigma}_0}{2}(1-u)^{\lambda-2}+o[(1-u)^{\lambda-2}] ,  &  0< \lambda <1 \\ \\
 \ds  {\tilde \sigma}_0\left(\frac{ \langle T \rangle_{\lambda} }{2} -1\right)(1-u)^{-1}+\frac{{\tilde \sigma}_0(1-u)^{\lambda-2}}{2H_{\lambda}(1)} 
+ o_I[(1-u)^{-1}] +o[(1-u)^{\lambda-2}] , &  \lambda >1  .
 \end{array} \right.
\eeq
The first line in this relation corresponds to standard Sibuya. In order to capture the leading contributions for large times, by Tauberian arguments, we have respectively picked up the lowest integer and non-integer orders in $1-u$
[Remark: For $1< \lambda <2$ since $-1 < \lambda -2 <0 $ we have $o_I(1-u)^{-1} = c_0+c_1(1-u)+\ldots \in o(1-u)^{\lambda-2}$ which is
consistent with our previous result \cite{SRW2022} -- see Eq. (26) and identify with $B_{\lambda} = -1/H_{\lambda}(1) >0$ for $\lambda \in (1,2)$].
This yields
\beq
\label{asymp_large-t}
\left\langle X_{\lambda}(t) \right\rangle \sim \left\{\begin{array}{lr} \ds \frac{\sigma_0}{2}\frac{(2-\lambda)_t}{t!} 
\to  \frac{{\tilde \sigma}_0}{2 \Gamma(2-\lambda)} t^{1-\lambda} , & \hspace{0.75cm} 0< \lambda < 1 \\ \\ \ds 
\frac{{\tilde \sigma}_0}{2}[\left\langle T\right\rangle_{\lambda} -2]  + {\tilde \sigma}_0\frac{(2-\lambda)_t}{2H_{\lambda}(1) t!} \to 
\frac{{\tilde \sigma}_0}{2}[\left\langle T\right\rangle_{\lambda} -2] + \frac{{\tilde \sigma}_0 t^{1-\lambda}}{2H_{\lambda}(1) \Gamma(2-\lambda)} 
& \hspace{0.75cm} \lambda > 1
\end{array}\right.     
\eeq 
For $\lambda \in (0,1)$ (standard Sibuya) the squirrel escapes to infinity along the direction of $\sigma_0$
by a $t^{1-\lambda}$-power law. 
Physically this can be interpreted by the occurrence of very long waiting times between the step reversals. 
\begin{figure}[t]
\centerline{\includegraphics[width=0.75\textwidth]{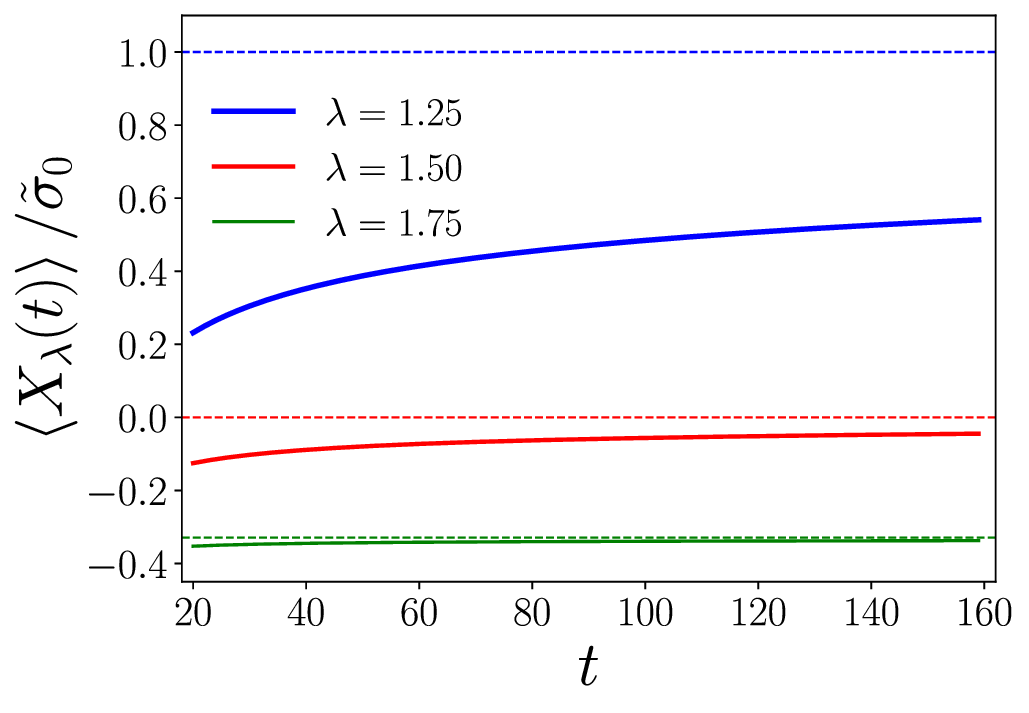}}
\vspace{-3mm}
\caption{\label{Fig2} Large time behavior of the expected position 
(\ref{final-result}) for $\lambda< \frac{3}{2}$, $\lambda=\frac{3}{2}$, and $\lambda>\frac{3}{2}$. The case $\lambda=\frac{3}{2}$ (red curve) corresponds to an asymptotically unbiased walk with $\langle X_{\frac{3}{2}}(\infty)\rangle = 0$. Dashed lines denote $\langle X_{\lambda}(\infty)\rangle$.}
\end{figure}
For narrower GSDs with $\lambda >1$ (shorter waiting times with existing mean $\langle T \rangle_{\lambda}$) the squirrel
remains trapped close to the departure site where the value 
$\langle X_{\lambda}(\infty)\rangle = \frac{{\tilde \sigma}_0}{2}( \langle T \rangle_{\lambda}  -2)$ is approached 
by a $t^{-(\lambda-1)}$-power law term which has opposite sign to ${\tilde \sigma}_0$ (see (\ref{final-result})). For $\lambda = \frac{3}{2}$ we have $\langle T\rangle_{1.5}=2$
and $\langle X_{1.5}(\infty)\rangle= 0$ where the walk is in the large time limit unbiased (in the average any second step is reversed). For $\lambda < \frac{3}{2} $ (i.e.\ $m=2$ and $\mu < 0.5$ with $\langle T\rangle_{\lambda} > 2$) the waiting times between the step reversals are still relatively long (the GSD being relatively broad) where $\langle X_{\lambda}(\infty)\rangle$
has the same sign as ${\tilde \sigma}_0$. In this case the squirrel does not escape in ${\tilde \sigma}_0$-direction, but in the average remains trapped on the same side of the departure site ($sign(\langle X_{\lambda}(\infty)\rangle) = {\tilde \sigma}_0$). 
This behavior changes for $\lambda > \frac{3}{2}$ which means shorter waiting times between the step switches and 
narrower GSD ($\langle T \rangle_{\lambda} < 2$): 
the sign of $\langle X_{\lambda}(\infty)\rangle$ changes and becomes opposite to ${\sigma}_0$. 
We can see this more closely if we rewrite (\ref{asymp_large-t}) for $\lambda >1$ ($\lambda=m-1 + \mu$ and $\mu\in (0,1)$) as
\beq
\label{final-result}
\begin{array}{cl} 
\ds \left\langle X_{\lambda}(t) \right\rangle &  \sim \ds  -\frac{{\tilde \sigma}_0}{2(m-2+\mu)}\left(m-3+2\mu  + \frac{\Gamma(t+3-m-\mu)}{\Gamma(t+1)}\frac{\Gamma(m)}{\Gamma(1-\mu)}\right) \\ \\
 & \to \ds  -\frac{{\tilde \sigma}_0}{2(m-2+\mu)}\left(m-3+2\mu +
\frac{\Gamma(m)}{\Gamma(1-\mu)}t^{-(m-2+\mu)}\right) 
 \end{array} \hspace{0.3cm} m=\ceil{\lambda} \geq 2
\eeq
This relation is plotted in Figure \ref{Fig2} for three values of $\lambda \in (1,2)$ including the asymptotically unbiased case $\lambda=\frac{3}{2}$ where the squirrel approaches the departure site with a
$-{\tilde \sigma}_0(\pi t)^{-\frac{1}{2}}$--law.

\section{Anomalous diffusive features}
\label{anomalous_transport}
In this section we analyze the mean square displacement (MSD) which we denote with
$\left\langle X_{\lambda}^2(t)\right\rangle $ (with respect to the initial position $X_{\lambda}(0)=0$) and especially focus on the large time asymptotics. The MSD is given by
\beq
\label{SRW_MSD}
\begin{array}{clr}
\ds \left\langle X_{\lambda}^2(t)\right\rangle &  = \ds \left\langle \sum_{r_1=1}^t\sum_{r_2=1}^t\sigma_{r_1}\sigma_{r_2} \right\rangle = \ds  2\sum_{r=1}^t \sum_{s=1}^r \left\langle \sigma_r\sigma_s  \right\rangle -\sum_{r=1}^t \left\langle \sigma_r^2  \right\rangle & \\ \\
& = \ds  2 K_{\lambda}(t) -t = - t + 2\sum_{r=1}^t\sum_{s=r}^t \left\langle (-1)^{{\cal N}_{\lambda}(r;s-r)} \right\rangle & \\ \\
& = \ds -t + 
2\sum_{r=0}^t\sum_{k=0}^{t-r}\left\langle (-1)^{{\cal N}_{\lambda}(r;k)} \right\rangle
-2\sum_{k=0}^{t}\left\langle (-1)^{{\cal N}_{\lambda}(0;k)} \right\rangle 
\end{array}  \hspace{1cm} t \in \mathbb{N}
\eeq
where comes into play the new quantity
\beq
\label{new_counting_var}
{\cal N}_{\lambda}(r;s-r) =
{\cal N}_{\lambda}(s)-{\cal N}_{\lambda}(r)  , 
\hspace{1cm} s \geq r \geq 0
\eeq
of the so called aged (generalized Sibuya) counting process ${\cal N}_{\lambda}(r;k)$ and ${\cal N}_{\lambda}(0;k)$ recovers the original GSP.
Aged renewal processes have been introduced and analyzed for continuous-time renewal processes \cite{BarkaiCheng2003,GodrecheLuck2001,SchulzBarkaiMetzler2014} and only recently for discrete time counting processes \cite{SRW2022}.
We emphasize that the aged renewal process 
${\cal N}(r;k)$ (apart of the Markovian cases, Bernoulli and Poisson) depends on the `aging parameter' $r$ and is different from the original 
counting process ${\cal N}(k)$ reflecting non-markovianity of the latter. 
To evaluate (\ref{SRW_MSD}) (see \cite{SRW2022} for more details)
it is useful to consider first the GF\footnote{We suppress here $\lambda$ in ${\cal N}_{\lambda}(t)$ to emphasize that this deduction holds for any discrete-time renewal process.}:
\beq
\label{second_term}
{\bar g}_v(w,u) = \sum_{\tau=0}^{\infty}\sum_{t=0}^{\infty}w^{\tau}u^t\left\langle v^{{\cal N}(\tau;t)}\right\rangle ,\hspace{1cm} |u|, |w| < 1,\hspace{0.5cm} |v|\leq 1.
\eeq
We further introduce the auxiliary function $h_v(r,t)=\Theta(t-r)\sum_{k=0}^{t-r}\left\langle v^{{\cal N}(r;k)} \right\rangle$ and its double GF
\beq
\label{Hv}
\begin{array}{clr}
\ds {\bar h}_v(w,u) &=  \ds \sum_{r=0}^{\infty}\sum_{t=0}^{\infty}w^{r}u^t 
\Theta(t-r)\sum_{k=0}^{t-r} \left\langle v^{{\cal N}(r;k)} \right\rangle  &\\ \\
& = \ds \sum_{s=0}^{\infty}u^s \sum_{k=0}^s \sum_{r=0}^{\infty}(wu)^{r} \left\langle v^{{\cal N}(r;k)}  \right\rangle      & \\ \\
 & = \ds \frac{{\bar g}_v(uw,u)}{1-u} &
\end{array}
\eeq
where in the second line we substitute $t=r+s$ and introduced the discrete Heaviside step function $\Theta(s) = 1$ for $s\geq 0$ and $\Theta(s)=0$ else (especially $\Theta(0)=1$).
We then can write for the MSD GF
\beq
\label{MSD_GF}
{\bar X}^{(2)}(u) = 2{\bar K}(u)-\frac{u}{(1-u)^2}
\eeq
where $\frac{u}{(1-u)^2} = \sum_{t=1}^{\infty} t u^t$ and with (\ref{Hv}) we have
\beq
\label{K_GF}
{\bar K}(u)=  {\bar h}_{-1}(1,u)- {\bar h}_{-1}(0,u) =  \frac{{\bar g}_{-1}(u,u)}{1-u}- \frac{{\bar g}_{-1}(0,u)}{1-u}.
\eeq
To evaluate this relation we need to determine ${\bar g}_v(w,u)$ which is the GF of the state polynomial of the aged counting process
\beq
\label{aged_process}
{\bar g}_v(w,u) = \sum_{t=0}^{\infty}\sum_{\tau=0}^{\infty}u^t w ^{\tau}\sum_{m=0}^{\infty}\mathbb{P}[{\cal N}(\tau;t)=m]v^m = \sum_{m=0}^{\infty}v^m {\bar \phi}^{(m)}_w(u)
\eeq
where ${\bar \phi}^{(m)}_w(u)$ stand for the double GF of the state probabilities $\mathbb{P}[{\cal N}(\tau;t)=m]$
of the aged process ${\cal N}(\tau;t)$ which we determine as
\beq
\label{double_GF_aged}
 \ds {\bar \phi}^{(m)}_w(u) =  \ds \sum_{t=0}^{\infty}\sum_{\tau=0}^{\infty} u^tw^{\tau} \left\{\begin{array}{ll} \ds \sum_{n=0}^{\infty}\left\langle\  
 \Theta(J_n,\tau,J_{n+1})\Theta(J_{n+m},t+\tau,J_{n+m+1})
 \right\rangle  , & m>0 \\   \\
 \ds \sum_{n=0}^{\infty}\left\langle\  
 \Theta(J_n,\tau,J_{n+1})\Theta(J_{n+1}-t-\tau-1)\right\rangle ,   & m=0.
\end{array}\right.
\eeq
Note that $\Theta(J_{n+1}-1 -t-\tau)$ indicates that $J_{n+1}-1\geq t+\tau$, i.e. that a state $n$ at time $\tau$ still persists at time $t+\tau$.
Using the IID feature of the $\Delta t$ and $\Theta(J_n,t+\tau,J_{n+1})= \Theta(J_n-\tau,t,J_{n+1}-\tau)$ this yields
\beq
\label{double_GF_aged_evaluation}
 {\bar \phi}^{(m)}_w(u) =  \left\{\begin{array}{lr} \ds \frac{u[{\bar \psi}(u)-{\bar \psi}(w)]}{(u-w)[1- {\bar \psi}(w)]}
{\bar \psi}(u)^{m-1}\frac{1-{\bar \psi}(u)}{1-u}  , & m>0 \\   \\
 \ds \frac{1}{(1-u)}
\left[\frac{1}{1-w}
-\frac{u[{\bar \psi}(u)-{\bar \psi}(w)]}{(u-w)[1-{\bar \psi}(w)]}\right]   ,   & m=0.
\end{array}\right.
\eeq
Now we can evaluate (\ref{aged_process}) to arrive at
\beq
\label{tripple-GF_result}
{\bar g}_v(w,u)=\left\{\begin{array}{lr} \ds   \frac{1}{(1-w)(1-u)} - \frac{(1-v)u}{(1-u)(u-w)[1-v{\bar \psi}(u)]}\frac{[{\bar \psi}(u)-{\bar \psi}(w)]}{[1-{\bar \psi}(w)]} , & u \neq w \\ \\ \ds  \frac{1}{(1-u)^2}
- \frac{(1-v)u}{(1-u)[1-v{\bar \psi}(u)][1-{\bar \psi}(u)]} \frac{d{\bar \psi}(u)}{du} , & u=w .
   \end{array}\right.
\eeq
This relation contains the GF of the state polynomial (\ref{state_gen})
${\bar g}_v(0,u) = \frac{1-{\bar \psi}(u)}{1-u}\frac{1}{1-v{\bar \psi}(u)}$
of the original counting process ${\cal N}(t)={\cal N}(0;t)$.
Then we obtain for (\ref{K_GF}) which determines the MSD GF (\ref{MSD_GF}) the expression
\beq
\label{F_genFuc}
{\bar K}(u) =  \frac{1}{(1-u)^3} -\frac{1}{(1-u)^2(1-[{\bar \psi}(u)]^2)}\left(2u\frac{d{\bar \psi}(u)}{du}+[1-{\bar \psi}(u)]^2\right).
\eeq
Now consider the large-time asymptotics of the MSD for the case when ${\cal N}(t)={\cal N}_{\lambda}(t)$ is the GSP
where we denote then (\ref{F_genFuc}) with ${\bar K}_{\lambda}(u)$. Using (\ref{utooneexpand}) for $u\to 1$ and
$ \frac{d{\bar \psi}_{\lambda}(u)}{du} = \sum_{t=1}^{\infty} t \psi_{\lambda}(t) u^{t-1}$ we have 
\beq
\label{utooneexpand_derivative}
\frac{d{\bar \psi}_{\lambda}(u)}{du} = \left\{\begin{array}{lr} \ds \lambda(1-u)^{\lambda-1} , & \hspace{0.3cm} 0 < \lambda < 1  \\ \\ \ds 
\left\langle T \right\rangle_{\lambda} +\frac{\lambda}{H_{\lambda}(1)}(1-u)^{\lambda-1} 
+ o_I(1)+o[(1-u)^{\lambda-1}] , & \hspace{0.3cm}  1 < \lambda < 2 \\ \\
 \ds   \left\langle T \right\rangle_{\lambda} -2 B_2^{(\lambda)}(1-u) +\frac{\lambda}{H_{\lambda}(1)}(1-u)^{\lambda-1} 
+ o_I(1) + o[(1-u)^{\lambda-1}] , & \lambda >2
\end{array}\right. (u \to 1-)
\eeq
where $o_I(1)$ is an expansion containing only integer powers $(1-u)^n$ with $n\geq 1$.
Thus we obtain ($u \to 1-$)
\beq
\label{K-lambda_GF}
\begin{array}{llr} 
\ds {\bar K}_{\lambda}(u) & = \ds  \frac{1}{(1-u)^3} 
-\frac{2 u}{(1-u)^2(1-{\bar \psi}_{\lambda}(u))(1+{\bar \psi}_{\lambda}(u))}\frac{d{\bar \psi}_{\lambda}(u)}{du}
- \frac{1}{(1-u)^2}\frac{1-{\bar \psi}_{\lambda}(u)}{1+{\bar \psi}_{\lambda}(u)} & \\ \\
 & = \ds   \left\{\begin{array}{lr} \ds 
 \frac{1}{(1-u)^3}\left(1-\frac{\lambda}{[1-\frac{1}{2}(1-u)^{\lambda}]}+\frac{\lambda(1-u)}{[1-\frac{1}{2}(1-u)^{\lambda}]}\right) -\frac{1}{2}\frac{(1-u)^{\lambda-2}}{[1-\frac{1}{2}(1-u)^{\lambda}]} \\ \\  \qquad = (1-\lambda)(1-u)^{-3}  +o[(1-u)^{-3}] , & \lambda \in (0,1)   \\ \\  \ds 
\left(1-\frac{\langle T\rangle_{\lambda}}{2}\right)\,(1-u)^{-2}  - \frac{\lambda-1}{H_{\lambda}(1) \langle T\rangle_{\lambda}}(1-u)^{\lambda-4} +o[(1-u)^{-2}] , & \hspace{1cm} \lambda \in (1,2) \\ \\
\ds
 (1-u)^{-2}\left(1+\frac{B_2^{(\lambda)}}{\langle T\rangle_{\lambda}} -\frac{ \langle T\rangle_{\lambda}}{2} \right) +o[(1-u)^{-2}] , & \hspace{1cm} \lambda > 2 
  \end{array}\right. & 
 \end{array} 
\eeq
where only the first two terms in (\ref{F_genFuc}) contain the relevant orders.
Inversion of (\ref{K-lambda_GF}) yields ($t\to \infty$)
\beq
\label{large_time_K}
K_{\lambda}(t) \sim  \left\{\begin{array}{lr} \ds \frac{(1-\lambda)}{2}t^2,  &  0<\lambda < 1 \\ \\ \ds 
\left(1-\frac{\langle T\rangle_{\lambda}}{2}\right)\, t- \frac{\lambda-1}{H_{\lambda}(1) \langle T\rangle_{\lambda}\Gamma(4-\lambda)} t^{3-\lambda} \to -  \frac{\lambda-1}{H_{\lambda}(1) \langle T\rangle_{\lambda}\Gamma(4-\lambda)} t^{3-\lambda}  , & 1<\lambda <2 \\ \\ \ds  \left(1+\frac{B_2^{(\lambda)}}{\langle T\rangle_{\lambda}} -\frac{ \langle T\rangle_{\lambda}}{2} \right)\,  t   , \hspace{1cm}
& \lambda >2
\end{array}\right.  
\eeq
The MSD (\ref{SRW_MSD}) then scales as
\beq
\label{MSD_largetime}
\langle X_{\lambda}^2(t)\rangle = 2K_{\lambda}(t)-t \sim  \left\{\begin{array}{lr} \ds (1-\lambda)t^2 \, , & 0<\lambda < 1 \\ \\
\ds  - 2\frac{\lambda-1}{H_{\lambda}(1) 
\langle T\rangle_{\lambda}\Gamma(4-\lambda)} t^{3-\lambda} \,  , & 1<\lambda <2 \\ \\ \ds 
 \left(1+\frac{2B_2^{(\lambda)}}{\langle T\rangle_{\lambda}} - \langle T\rangle_{\lambda} \right)\, t = \frac{t}{\langle T \rangle_{\lambda}}\left(\langle T^2 \rangle_{\lambda} - (\langle T \rangle_{\lambda})^2\right) , & \lambda >2
\end{array}\right.
\eeq
where all quantities are non-negative and with the GSD variance ${\cal V}_{\lambda}=\langle T^2 \rangle_{\lambda} - (\langle T \rangle_{\lambda})^2$ determined in (\ref{variance_GSD}) and the mean waiting time $\langle T \rangle_{\lambda}$ in (\ref{first_moment}).
Hence this relation writes
\beq
\label{MSD_largetime_GSP}
\langle X_{\lambda}^2(t) \rangle \sim  \left\{\begin{array}{lr} \ds (1-\lambda)t^2 \, , & 0<\lambda < 1 \\ \\
\ds  \frac{2(\lambda-1)}{\Gamma(4-\lambda)} t^{3-\lambda} \,  , & 1<\lambda <2 \\ \\ \ds 
\frac{\lambda(m-\lambda)}{(\lambda-1)(\lambda-2)}  t \,  , & \lambda >2.
\end{array}\right. \hspace{0.5cm} (t \to \infty)
\eeq
In view of the power-laws governing the expected position (\ref{asymp_large-t}) one can see that 
$\langle X_{\lambda}^2(t) \rangle \gg \langle X_{\lambda}(t) \rangle^2$. Therefore, the MSD (\ref{MSD_largetime_GSP}) dominates the large-time asymptotics of the spatial variance of the squirrel motion.
The normal diffusive behavior occurring for $\lambda>2$ breaks down at the limits $\lambda=m-$ ($m\geq 3$) where we have 
${\cal V}_{m-} = 0$ (see (\ref{evalB}) -- (\ref{variance_GSD})) with deterministic oscillatory squirrel motions.
Contrarily to these cases 
the limit $\lambda=2-$ is non-deterministic which is expressed by $\langle X_{2-}^2(t) \rangle = 2t$ (see (\ref{MSD_largetime_GSP})) corresponding to persistent normal diffusion (Brownian motion) of the squirrel with spatial Gaussian limiting distribution of propagator (\ref{rewrite_SRW_propagator}).
The limiting cases $\lambda=n-$ and $\lambda=n+$ ($n\in \mathbb{N}$) exhibiting respectively distinct behaviors are considered in the Appendix more closely.
\\[1mm]
We identify three different diffusive large-time regimes for the {\it generalized Sibuya SRW}:
\begin{enumerate}
\item[(i)] A ballistic superdiffusive regime when the GSD is broad with $0< \lambda < 1$ (standard Sibuya) with a $t^2$--law. 
\item[(ii)] A superdiffusive regime for $1<\lambda<2$ with a $t^{3-\lambda}$--law with scaling exponent $ 1< 3-\lambda < 2$.
\item[(iii)] A normal diffusive regime when the GSD is narrow for $\lambda >2$ with emergence of Brownian motion.
\end{enumerate}
These results are consistent with those obtained in our recent paper \cite{SRW2022} by considering general asymptotic features of discrete-time renewal processes. 
Superdiffusive large time regimes of these types
were also reported for continuous time Cattaneo transport models \cite{CompteMetzler1997,CompteMetzler1999}.
The generalized Sibuya SRW of the present study covers for different ranges of $\lambda$ the whole spectrum from anomalous--ballistic (i), over anomalous (ii) to normal (Brownian) diffusion (iii).

\section{Conclusions}
\label{Conclusions}
In the present paper we have studied a semi-Markovian discrete-time generalization of the telegraph (Cattaneo) process
where the waiting times between the step reversals follow the {\it generalized Sibuya distribution} -- GSD. We called this walk the {\it generalized Sibuya SRW}. It turns out that the presented model has a large flexibility to cover a wide range of behaviors including superdiffusive-ballistic, superdiffusive and normal diffusive transport. We have shown that these features are solely governed by the ``broadness'' of the GSD waiting time density.
For follow-up research an interesting subject is the analysis of scaling limits to continuous time and space which 
define new semi-Markovian generalizations of telegraph (Cattaneo) processes.

Moreover, variants of SRW models in multidimensional spaces appear to be interesting directions. For instance in problems where a walker is moving with constant velocity in a $D$-dimensional infinite space and changing its velocity direction randomly at the renewal times of a discrete-time counting process such as the GSD or others. 
The class of generalized Sibuya SRW and similar models open various new directions in random walk theory, general fractional calculus and non-Markovian  dynamics in complex systems. 
\subsubsection*{Acknowledgments}

F.~Polito has been partially supported by INdAM/GNAMPA.
%
%
%
%
%
{
\appendix
\setcounter{equation}{0}
\renewcommand{\theequation}{A.\arabic{equation}}
\section{Appendix}
\label{Appendix}
Let us discuss here the behaviors emerging in the limiting cases $\lambda \to n-$ and $\lambda \to n+$ ($n\in \mathbb{N}$), respectively. \\
(a) $\lambda= 2+0$ ($\lambda=2+\epsilon$ with $\epsilon \to 0+$ and $m=3$):
Then we have $\langle T \rangle_{2+\epsilon} = \frac{2}{1+\epsilon} \sim 2- \neq \langle T \rangle_{2-\epsilon} = \frac{1}{1-\epsilon} \sim 1+$ and
${\cal V}_{2+\epsilon} \sim \frac{(2+\epsilon)(1-\epsilon)}{(1+\epsilon)^2}\frac{2}{\epsilon} \sim \frac{4}{\epsilon}$ thus
\beq
\label{var_two_plus}
\frac{{\cal V}_{2+\epsilon}}{\langle T\rangle_{2+\epsilon}} \sim \frac{2}{\epsilon} \to \infty
\eeq
i.e. for $\lambda=2+$ the MSD (\ref{MSD_largetime_GSP}) $$\langle X_{2+\epsilon}^2(t) \sim  \langle X_{2+\epsilon}^2(t) -\langle X_{2+\epsilon}(t)\rangle^2 \sim \rangle \sim \frac{2}{\epsilon} t \to \infty $$ is singular where the average position $\langle X_{2+\epsilon}(t)\rangle \sim -\frac{{\tilde \sigma}_0\epsilon}{1+\epsilon} \to 0$, see (\ref{final-result}). In the limit
$\lambda=2+$ emerges for large observation times Brownian diffusion (according to case (iii)) where the squirrel position is in the average on the departure site, but with extremely large fluctuations. We observe in (\ref{MSD_largetime_GSP}) that the normal diffusive behavior is not singular at the limits $\lambda=n+$ ($n\in \mathbb{N}$) for $n>2$ and it is also different at the limit $\lambda=2-$ which we consider next.
\\[2mm]
(b) $\lambda= 2-0$ ($\lambda=2-\epsilon$ with $\epsilon \to 0+$ and $m=2$):
Let us compare this limit with (a).
We then have $H_{2-\epsilon}(1) = 
\frac{\Gamma(\epsilon)}{\Gamma(\epsilon-1)\Gamma(2)}=\epsilon-1$ and $\mu=1-\epsilon$ with (\ref{MSD_largetime_GSP})
\beq
\label{MSD-2minus}
\langle X_{2-}^2 \rangle
\sim \langle X_{2-}^2 \rangle - \langle X_{2-} \rangle^2   \sim   2\frac{1-\epsilon}{\Gamma(2+\epsilon)} t^{1+\epsilon} \to 2t
\eeq
where with (\ref{final-result}) we see that $\langle X_{2-} \rangle^2 =\frac{1}{4} \ll \langle X_{2-\epsilon}^2 \rangle$.
This limit corresponds to Brownian motion of the squirrel and is different from the (deterministic) trivial oscillatory motion of the case $\lambda=2$ and is also different from the fast Brownian motion emerging in the limit $\lambda=2+$ of (a).
\\[2mm]
It is worthy of mention that relation (\ref{MSD-2minus}) can be re-derived in the following different way.
GF (\ref{gen-Sib}) has for $\lambda=2-$ the form
\beq
\label{GF2minus}
{\bar \psi}_{2-\epsilon}(u) =\frac{1}{(\epsilon-1)u}[\epsilon-1+(2-\epsilon)(1-u)-(1-u)^{2-\epsilon}]
\eeq
with small $\epsilon>0$. For $u\to 1-$ we have
\beq
\label{asymppsi}
{\bar \psi}_{2-\epsilon}(u) = 1 -\frac{1}{1-\epsilon}(1-u)+ \frac{1}{1-\epsilon}(1-u)^{2-\epsilon} +o(1-u).
\eeq
Now consider the second derivative
\beq
\label{consider}
{\bar b}_{\epsilon}(u) = \frac{1}{2} \frac{d^2}{du^2} {\bar \psi}_{2-\epsilon}(u) \sim \frac{1}{2} (2-\epsilon)(1-u)^{-\epsilon}  \sim(1-u)^{-\epsilon}
\eeq
inversion yields a Dirac $\delta$-distribution $b_{\epsilon}(t) =\frac{(\epsilon)_t}{t!}  \sim \frac{t^{\epsilon-1}}{\Gamma(\epsilon)} \to \delta_{+}(t)$ concentrated at $t=0+$.
Although the second derivative does not exist at $u=1$ we can define it in a distributional (Gel’fand-Shilov) sense
\cite{GelfandShilov1968}
to define constant $B_2^{(\lambda)}$ for the (`forbidden') limit 
\beq
\label{forbidden_limit}
B_2^{(2-)}= {\bar b}_{0+}(1-) \sim
  \int_{0}^{\infty} \delta_{+}(t){\rm d}t = 1.
\eeq
We hence have for the GSD variance
\beq
\label{GSD_var-2minus}
{\cal V}_{2-} \sim 2 B_2^{(2-)} + \langle T \rangle_{2-} -\langle T \rangle_{2-}^2 = 2 
\eeq
which takes us with (\ref{MSD_largetime}) and $\langle T \rangle_{2-} = 1+$ back to relation (\ref{MSD-2minus}). 
\\ 
In fact what we are using in (\ref{asymppsi}) is that the asymptotic expansion of ${\bar \psi}_{2-\epsilon}(u)$ for $u\to 1-$ captures the dominating contribution
of the GSD power-law tail ${\bar \psi}_{2-\epsilon}(t) \sim \frac{t^{\epsilon-3}}{(1-\epsilon)\Gamma(\epsilon-2)}$, namely
\beq
\label{generlized_sense}
{\cal V}_{2-} \sim 2B_2^{(2-)}   \sim  
\frac{1}{1-\epsilon)}\int_0^{\infty} \tau(\tau-1) \frac{\tau^{\epsilon-3}}{\Gamma(\epsilon-2)}  {\rm d}\tau 
\sim \int_0^{\infty}\delta_{+}(\tau)\frac{d^2}{d\tau^2}\tau^2 {\rm d}\tau =
 2
\eeq
where we use G`elfand-Shilov distributional relation $\frac{\tau^{\epsilon-3}}{\Gamma(\epsilon-2)} \to \frac{d^2}{d\tau^2}\delta_{+}(\tau)$ which only captures the information of the highest moment $\langle T^2 \rangle_{2-\epsilon}$. The contribution (\ref{generlized_sense}) is in a sense due to the power law tail of $ \frac{\tau^{\epsilon-3}}{\Gamma(\epsilon-2)}$ of $\psi_{2-\epsilon}(t)$
which is dying out for $\epsilon\to 0+$ and which is null for the deterministic case with the exact value $\lambda=2$.
\\[2mm]
In the same way we can consider the $m$th moment in the limit $\lambda=m-$ for any $m\in \mathbb{N}$.
We then have for the tail (\ref{large_time_pdf}) the distributional relation
\beq
\label{limit_minus}
\psi_{m-\epsilon}(t) \sim \frac{(m-\epsilon)[(m-1)!]}{{\Gamma(\epsilon)}} t^{\epsilon-1-m} \sim
(-1)^m \frac{d^m}{dt^m} \frac{t^{\epsilon-1}}{\Gamma(\epsilon)}
\to (-1)^m\frac{d^m}{dt^m}\delta_{+}(t)
\eeq 
which leads to the finite limiting value for the $m$th moment
\beq
\label{asympt_m_moments}
\langle T^m \rangle_{m-} \sim \int_0^{\infty} \tau^m (-1)^m\frac{d^m}{d\tau^m}\delta_{+}(\tau){\rm d}\tau = \Gamma(m+1) \, ,\hspace{1cm} m \in \mathbb{N}.
\eeq
For $m=1$ this yields $\langle T \rangle_{1-} = 1$ and is different from singular limiting case 
$\langle T \rangle_{1+\epsilon} = \frac{1}{\epsilon} \to \infty$. 
For a further discussion of the (standard Sibuya) limit $\lambda=1-$, we refer to our recent paper \cite{SRW2022}.
\\[2mm]
Now it is instructive to compare $\lambda=n-$ with
$\lambda=n+$ where $n\in \mathbb{N}$. In the latter case we have
$m=\ceil{n+}=n+1$ thus (\ref{large_time_pdf}) takes the form
\beq
\label{right-limit}
\psi_{n+\epsilon}(t) \sim \frac{(n+\epsilon)\Gamma(n+1)}{\Gamma(1-\epsilon)} t^{-n-1-\epsilon} 
\eeq
which remains `broad' behaving as $t^{-1-n}$ when $\epsilon \to 0+$ (contrarily to the limit $\lambda=n-$, see (\ref{limit_minus}) for $n=m$).
Therefore the $n$th moment 
\beq
\label{nth_moment}
\langle T^n \rangle_{n+\epsilon} \sim \int_0^{\infty} \tau^n \psi_{n+\epsilon}(\tau){\rm d}\tau \sim \frac{n \Gamma(n+1)}{\epsilon} \to \infty , \hspace{1cm} n \in \mathbb{N}
\eeq
has a $1/\epsilon$-singularity.
\\[2mm]
Considering now again $n=2$ the two limits $\lambda=2-$ and $\lambda=2+$, their difference becomes clear when we look at
the tails of the GSD (\ref{large_time_pdf}) where $m(2+)=3=m(2-)+1$.
Finally we have
\beq
\label{psi_two_plus}
{\bar \psi}_{2+\epsilon}(t) \sim \frac{2(2+\epsilon)}{\Gamma(1-\epsilon)} t^{-3-\epsilon} .
\eeq
Contrary to $\psi_{2-0}(t)$ the PDF (\ref{psi_two_plus}) remains broad for $\epsilon \to 0$. Therefore,
\beq
\label{meansquare}
\langle T^2 \rangle_{2+\epsilon} \sim \frac{2(2+\epsilon)}{\Gamma(1-\epsilon)} \int_0^{\infty}  \tau^{-3-\epsilon} \tau^2{\rm d}\tau = \frac{2(2+\epsilon)}{\Gamma(1-\epsilon)} \frac{\tau^{-\epsilon}}{(-\epsilon)}\big|_0^{\infty} \sim \frac{4}{\epsilon}.
\eeq
Then we further have $\langle T \rangle_{2+\epsilon}=\frac{2}{1+\epsilon} \to 2- \ll \langle T^2 \rangle_{2+\epsilon} $ 
thus
\beq
\label{asymp}
\frac{{\cal V}_{2+\epsilon}}{\langle T \rangle_{2+\epsilon}} 
\sim \frac{\langle T^2 \rangle_{2+\epsilon}}{\langle T \rangle_{2+\epsilon}} \sim \frac{2}{\epsilon}
\eeq
bringing us back to (\ref{var_two_plus}).
}
%

\end{document}